\documentclass[12pt]{article}
\usepackage{amsmath,amssymb}
\usepackage{amscd}
\usepackage[all]{xy}
\newtheorem{thm}{Theorem}
\newtheorem{df}{Def\/inition}
\newtheorem{prop}{Proposition}

\newtheorem{lemma}{Lemma}

\DeclareMathOperator{\ord}{ord} 
\DeclareMathOperator{\Lie}{Lie}

\DeclareMathOperator{\GL}{\pmb{GL}} 
\DeclareMathOperator{\depth}{depth} 
\DeclareMathOperator{\Image}{Im} 
 
\DeclareMathOperator{\id}{id}
 
\def\endproof{$\hfill \square$}

\begin{document}
\title{Breuil's classification of $p$-divisible groups over regular
  local rings of arbitrary dimension}
\author{Adrian Vasiu and Thomas Zink}
\maketitle

\medskip\noindent
To appear in {\it Advanced Studies in Pure Mathematics}, Algebraic and Arithmetic Structures of Moduli Spaces (Sapporo 2007)

\medskip\noindent
{\bf Abstract.} Let $k$ be a perfect field of characteristic $p \geq
3$. We classify  $p$-divisible groups
over regular local rings of the form
$W(k)[[t_1,\ldots,t_r,u]]/(u^e+pb_{e-1}u^{e-1}+\ldots+pb_1u+pb_0)$,
where $b_0,\ldots,b_{e-1}\in W(k)[[t_1,\ldots,t_r]]$ and $b_0$ is an
invertible element. This classification was in the case $r = 0$
conjectured by Breuil and proved by Kisin.   


\section{Introduction}
Let $p\in\mathbb{N}$ be an odd prime. Let $k$ be a perfect field of
characteristic $p$. Let $W(k)$ be the ring of Witt vectors with
coefficients in $k$. Let $r\in\mathbb{N}\cup\{0\}$.  We consider the
ring of formal power series
\begin{displaymath}
  \mathfrak{S}:=W(k)[[t_1,\ldots,t_r,u]].
\end{displaymath}
We extend the Frobenius endomorphism $\sigma$ of $W(k)$ to
$\mathfrak{S}$ by the rules
\begin{equation}\label{VZ0e-1}
  \sigma(t_i)=t^p_i \quad \mbox{ and} \quad \sigma(u)=u^p.
\end{equation}
If $M$ is a $\mathfrak{S}$-module we define

\begin{displaymath}
  M^{(\sigma)} := \mathfrak{S} \otimes_{\sigma, \mathfrak{S}} M.  
\end{displaymath}
Let $e\in\mathbb{N}$. Let
\begin{displaymath}
  E(u)=u^e+a_{e-1}u^{e-1}+\dots+ a_1 u+a_0
\end{displaymath}
be a polynomial with coefficients in $W(k)[[t_1,\dots, t_r]]$ such
that $p$ divides $a_i$ for all $i\in\{0,\ldots,e-1\}$ and moreover
$a_0/p$ is a unit in $W(k)[[t_1,\dots, t_r]]$.

We define
\begin{displaymath}
  R := \mathfrak{S}/E\mathfrak{S};
\end{displaymath}
it is a regular local ring of dimension $r+1$ with parameter system
$t_1,\dots, t_r, u$.

The following notion was introduced in [B] for the case $r = 0$.

\begin{df} A Breuil window relative to $\mathfrak{S} \rightarrow R$ is
  a pair $(Q,\phi)$, where $Q$ is a free $\mathfrak{S}$-module of
  finite rank and where $\phi: Q\to Q^{(\sigma)}$ is a
  $\mathfrak{S}$-linear map whose cokernel is annihilated by $E$.
\end{df}
As $\phi[\frac{1}{E}]: Q[\frac{1}{E}]\to Q^{(\sigma)}[\frac{1}{E}]$ is
a $\mathfrak{S}[\frac{1}{E}]$-linear epimorphism between free
$\mathfrak{S}[\frac{1}{E}]$-modules of the same finite rank, it is an
injection. This implies that $\phi$ itself is an injection. We check
that $C:=\text{Coker}(\phi)$ is a free $R$-module. For this, we can
assume that $C\neq 0$ and thus that the $\mathfrak{S}$-module $C$ has
projective dimension $1$. Since $\mathfrak{S}$ is a regular ring, we
have $\depth C =\dim \mathfrak{S} - 1 = \dim R$. As $C$ has the same
depth viewed as an $R$-module or as a $\mathfrak{S}$-module, we
conclude that $C$ is a free $R$-module.

The goal of the paper is to prove the following result whose validity
is suggested by previous works of Breuil and Kisin (see [B] and [K]).

\begin{thm}\label{VZ1t} The category of $p$-divisible groups over $R$ is
  equivalent to the category of Breuil windows relative to
  $\mathfrak{S} \rightarrow R$.
\end{thm}
This theorem was proved by Kisin [K] in the case $r=0$. We prove the
generalization by a new method which is based on the theory of
Dieudonn\'e displays \cite{ZDD}. This theory works only for a perfect
field $k$ of characteristic $p \geq 3$. 

Our method yields results for more general fields if we restrict
ourselves to formal 
$p$-divisible groups over $R$, i.e. $p$-divisible groups over $R$ whose special
fibers over $k$ are connected. Then we can substitute Dieudonn\'e displays by
nilpotent displays. To state the result we have to define nilpotent
Breuil windows relative to $\mathfrak{S} \rightarrow R$. Let $(Q, \phi)$ be a Breuil window relative to $\mathfrak{S} \rightarrow R$. Then we define a
$\sigma$-linear map $F: Q^{(\sigma)} \rightarrow Q^{(\sigma)}$ by
$F(x) = \id \otimes \phi^{-1}(Ex)$ for $x \in Q^{(\sigma)}$. If we tensor
$(Q^{(\sigma)}, F)$ by the $W(k)$-epimorphism $\mathfrak{S} \rightarrow W(k)$
which maps the variables $t_i$ and $u$ to $0$, we obtain a Dieudonn\'e
module over $W(k)$. In the theorem above this is the covariant
Dieudonn\'e module of the special fibre of the $p$-divisible group over $R$ which
corresponds to $(Q,\sigma)$. We say that $(Q, \phi)$ is a nilpotent Breuil window relative to $\mathfrak{S} \rightarrow R$ if
$(Q^{(\sigma)}, F) \otimes_{\mathfrak{S}} W(k)$ is the covariant
Dieudonn\'e module of a connected $p$-divisible group over $k$. Then
the arguments of this paper show that for $p = 2$ the category of
nilpotent Breuil windows is equivalent to the category of formal
$p$-divisible groups over $R$. Exactly the same statement holds for
non perfect fields of arbitrary characteristics $p>0$ if we replace
the ring $W(k)$ by a Cohen ring $C_k$ and $\mathfrak{S}$ by
$C_k[[t_1,\ldots t_r,u]]$. 

One can view Theorem \ref{VZ1t} as a ramified analogue of Faltings
deformation theory over rings of 
the form $W(k)[[t_1,\ldots,t_r]]$ (see [F, Thm. 10]). The importance
of Theorem 1 stems from its potential applications to modular and
moduli properties and aspects of Shimura varieties of Hodge type (see
\cite{VZ} for applications with $r=1$).

As in the case $r=0$ (see \cite{K}), Theorem 1 implies a
classification of finite flat, commutative group schemes of $p$ power
order over $R$.
\begin{df}
  A Breuil module relative to $\mathfrak{S} \rightarrow R$ is a pair
  $(M, \varphi)$, where $M$ is a $\mathfrak{S}$-module of projective
  dimension at most one and annihilated by a power of $p$ and where
  $\varphi:M \rightarrow M^{(\sigma)}$ is a $\mathfrak{S}$-linear map
  whose cokernel is annihilated by $E$.
\end{df}

\begin{thm}\label{VZ2t}
  The category of finite flat, commutative group schemes of $p$ power
  order over $R$ is equivalent to the category of Breuil modules
  relative to $\mathfrak{S}\to R$.
\end{thm}

The first author would like to thank MPI Bonn, Binghamton University,
and Bielefeld University for good conditions to work on the paper. The
second author would like to thank Eike Lau for helpful discussions.
Both authors thank the referee for some valuable remarks. 

\section{Breuil windows modulo powers of $u$}

We need a slight variant of Breuil windows, which was also considered
by Kisin in his proof of Theorem \ref{VZ1t} for $r=0$.

For $a\in \mathbb{N}$ we define $\mathfrak{S}_a :=
\mathfrak{S}/(u^{ae})$; it is a $p$-adic ring without $p$-torsion.
Clearly $E$ is not a zero divisor in $\mathfrak{S}_a$. The Frobenius
endomorphism $\sigma $ of $\mathfrak{S}$ induces naturally a Frobenius
endomorphism $\sigma $ of $\mathfrak{S}_a$.

We write
\begin{equation}\label{VZ0e}
  E = E(u) = u^e + p \epsilon,
\end{equation} 
where $\epsilon:=(a_{e-1}/p)u^{e-1} + \ldots + (a_1/p) u + (a_0/p)$ is
a unit in $\mathfrak{S}$. As $u^{ae}$ and $p^a(-\epsilon)^a$ are
congruent modulo the ideal $(E)$, we have identities
\begin{displaymath}
  \mathfrak{S}_a/(E) = \mathfrak{S}/(E,p^a) = R/p^aR.
\end{displaymath}

\begin{df} A Breuil window relative to $\mathfrak{S}_a\to R/p^aR$ is a pair $(Q,\phi)$, where $Q$ is a free
  $\mathfrak{S}_a$-module of finite rank and where $\phi: Q\to
  Q^{(\sigma)}$ is a $\mathfrak{S}$-linear map whose cokernel is
  annihilated by $E$ and is a free $R/p^aR$-module.
\end{df}
We will call this shortly a $\mathfrak{S}_a$-window, even though this
is not a window in the sense of \cite{ZW}. To avoid cases,
we define $\mathfrak{S}_{\infty}:= \mathfrak{S}$ and
$R/p^{\infty}R:=R$ and we will allow $a$ to be $\infty$. Thus from now
on $a\in \mathbb{N} \cup\{\infty\}$. We note down that a
$\mathfrak{S}_{\infty}$-window will be a Breuil window relative to
$\mathfrak{S}\rightarrow R$.  Next we relate $\mathfrak{S}_a$-windows
to the windows introduced in [Z3].

We will use the following convention from \cite{ZDFG}. Let $\alpha : M
\rightarrow N$ be a $\sigma$-linear map of $\mathfrak{S}_a$-modules. Then
we denote by 
$$\alpha^{\sharp} : M^{(\sigma)} = \mathfrak{S}_a \otimes_{\sigma,
  \mathfrak{S}_a} M \rightarrow N $$
its linearization. We say that $\alpha$ is a $\sigma$-linear
epimorphism (etc.) if $\alpha^{\sharp}$ is an epimorphism.

We consider triples of the form $(P,Q,F)$, where $P$ is a free
$\mathfrak{S}_a$-module of finite rank, $Q$ is a
$\mathfrak{S}_a$-submodule of $P$, and $F:P\to P$ is a $\sigma$-linear
map, such that the following two properties hold:

\begin{enumerate}
\item[(i)] $E\cdot P\subset Q$ and $P/Q$ is a free $R/p^aR$-module.

  \smallskip
\item[(ii)] $F(Q)\subset \sigma(E)\cdot P$ and $F(Q)$ generates
  $\sigma(E)\cdot P$ as a $\mathfrak{S}_a$-module.
\end{enumerate}
As $\sigma(E)$ is not a zero divisor in $\mathfrak{S}_a$, we can
define $F_1 := (1/\sigma(E))F : Q \rightarrow P$.

Any triple $(P,Q,F)$ has a \emph{normal decomposition}. This means
that there exist $\mathfrak{S}_a$-submodules $J$ and $L$ of $P$ such
that we have:
\begin{equation}\label{VZ-1e}
  P = J \oplus L\;\;\mbox{and}\;\;\quad Q = E\cdot J \oplus L. \end{equation}
This decomposition shows that $Q$ is a free $\mathfrak{S}_a$-module.
The map 
\begin{equation}\label{VZ0e-1a}
  F \oplus F_1: J \oplus L \rightarrow P
\end{equation} 
is a $\sigma$-linear isomorphism. A normal decomposition of $(P,Q,F)$
is not unique.

If $\tilde{P}$ is a free $\mathfrak{S}_a$-module of finite rank and if
$\tilde{P} = \tilde{L} \oplus \tilde{J}$ is a direct sum
decomposition, then each arbitrary $\sigma$-linear isomorphism
$\tilde{J} \oplus \tilde{L} \rightarrow \tilde{P}$ defines naturally a
triple $(\tilde{P},\tilde{Q},\tilde{F})$ as above. We can often
identify the triple $(\tilde{P},\tilde{Q},\tilde{F})$ with an
invertible matrix with coefficients in $\mathfrak{S}_a$ which is a
matrix representation of the $\sigma$-linear isomorphism $\tilde{J}
\oplus \tilde{L} \rightarrow \tilde{P}$. Each triple $(P,Q,F)$ is
isomorphic to a triple constructed as
$(\tilde{P},\tilde{Q},\tilde{F})$.

\begin{lemma}\label{Tripel}
  The category of triples $(P,Q,F)$ as above is equivalent to the
  category of $\mathfrak{S}_a$-windows.
\end{lemma} 
{\bf Proof:} Assume we are given a triple $(P,Q,F)$. By
definition $F_1$ induces a $\mathfrak{S}_a$-linear epimorphism
\begin{equation}\label{VZ0e-1aa}
  F_1^{\sharp}: Q^{(\sigma)}=\mathfrak{S}_a \otimes_{\sigma,
    \mathfrak{S}_a} Q\twoheadrightarrow P.   
\end{equation}
Due to the existence of normal decompositions of $(P,Q,F)$, $Q$ is a
free $\mathfrak{S}_a$-module of the same rank as $P$. Therefore $
F_1^{\sharp}$ is in fact an isomorphism. To the triple $(P,Q,F)$ we
associate the $\mathfrak{S}_a$-window $(Q,\phi)$, where
\begin{displaymath}
  \phi : Q \rightarrow Q^{(\sigma)}
\end{displaymath} 
is the composite of the inclusion $Q \subset P$ with $
(F_1^{\sharp})^{-1}$.

Conversely assume that we are given a $\mathfrak{S}_a$-window
$(Q,\phi)$.  We set $P := Q^{(\sigma)}$ and we consider $Q$ as a
submodule of $P$ via $\phi$. We denote by $F_1: Q \rightarrow P$ the
$\sigma $-linear map which induces the identity $Q^{(\sigma)} =
P$. Finally we set $F (x) := F_1 (Ex)$ for $x \in P$. Then $(P,Q,F)$
is a triple as above. \endproof

\smallskip Henceforth we will not distinguish between triples and
$\mathfrak{S}_a$-windows i.e., we will identify $(Q,\phi)\equiv
(P,Q,F)$. We can describe a normal decomposition directly in terms of
$(Q,\phi)$. Indeed, we can identify $Q$ with $J \oplus L$ via
$(1/E)\id_J \oplus \id_L$. Then a normal decomposition of $(Q,\phi)$
is a direct sum decomposition $Q = J \oplus L$ which induces a normal
decomposition $P=Q^{(\sigma)} = (1/E)\phi(J) \oplus \phi(L)$. If $a\in
\mathbb{N}$, then each $\mathfrak{S}_a$-window lifts to a
$\mathfrak{S}_{a+1}$-window (this is so as each invertible matrix with
coefficients in $\mathfrak{S}_a$ lifts to an invertible matrix with
coefficients in $\mathfrak{S}_{a+1}$).

\section{The $p$-divisible group of a Breuil window}\label{Sect3}

We relate $\mathfrak{S}_a$-windows to Dieudonn\'e displays over
$R/p^aR$ as defined in \cite{ZDD}, Definition 1. Let $S$ be a complete
local ring with residue field $k$ and maximal ideal $\mathfrak{n}$.
We denote by $\hat{W}(\mathfrak{n})$ the subring of all Witt vectors
in $W(\mathfrak{n})$ whose components converge to zero in the
$\mathfrak{n}$-adic topology. From \cite{ZDD}, Lemma 2 we get that
there exists a unique subring $\hat{W}(S) \subset W(S)$, which is
invariant under the Frobenius $F$ and Verschiebung $V$ endomorphisms
of $W(S)$ and which sits in a short exact sequence:
\begin{displaymath}
  0 \rightarrow \hat{W}(\mathfrak{n}) \rightarrow \hat{W}(S)
  \rightarrow W(k) \rightarrow 0.
\end{displaymath}
It is shown in \cite{ZDD} that the category of $p$-divisible groups
over $S$ is equivalent to the category of Dieudonn\'e displays over
$\hat{W}(S)$.

For $a\in \mathbb{N} \cup\{\infty\}$ there exists a unique
homomorphism
\begin{equation}
  \delta_a : \mathfrak{S}_a \rightarrow \hat{W}(\mathfrak{S}_a)
\end{equation}
such that for all $x \in \mathfrak{S}_a$ and for all $n\in \mathbb{N}$
we have $\mathbf{w}_n(\delta_a(x)) = \sigma^n(x)$ (here $\mathbf{w}_n$
is the $n$-th Witt polynomial). It maps $t_i \mapsto
[t_i]=(t_i,0,0,\ldots)$ and $u \mapsto [u]=(u,0,0,\ldots)$.  If we
compose $\delta_a$ with the canonical $W(k)$-homomorphism
$\hat{W}(\mathfrak{S}_a) \rightarrow \hat{W}(R/p^aR)$ we obtain a
$W(k)$-homomorphism
\begin{equation}\label{VZ0e-1b}
  \begin{array}{lccr}
    \varkappa_a: & \mathfrak{S}_a & \rightarrow & \hat{W}(R/p^aR).
  \end{array}
\end{equation}
We note that $p$ is not a zero divisor in $\hat{W}(R)$.
\begin{lemma}
  The element $\varkappa_{\infty}(\sigma(E)) \in \hat{W}(R)$ is
  divisible by $p$ and the fraction \\[2mm] 
$\tau :=  \varkappa_{\infty}(\sigma(E))/p$ is a unit in $\hat{W}(R)$.
\end{lemma} {\bf Proof:} We have $\varkappa_{\infty}(E) \in
V\hat{W}(R)$. Since $\varkappa_{\infty}$ is equivariant with respect
to $\sigma$ and the Frobenius $F$ of $\hat{W}(R)$ we get:
\begin{displaymath}
  \varkappa_{\infty}(\sigma(E)) = F(\varkappa_{\infty}(E)) \in p\hat{W}(R).
\end{displaymath}
We have to verify that $\mathbf{w}_0(\tau)$ is a unit in $R$. We have:
\begin{displaymath}
  \mathbf{w}_0(\tau) = \mathbf{w}_0( \varkappa_{\infty}(\sigma(E)))/p =
  \sigma(E)/p.
\end{displaymath}
With the notation of (\ref{VZ0e}) we have $\sigma (E) = u^{ep} + p\sigma
(\epsilon) $. Since $u^{ep} \equiv (p\epsilon)^p \; mod\; (E)$ we see
that $\sigma(E)/p$ is a unit in $R$. We note that this proof works for all primes $p$ (i.e., even if $p=2$). \endproof

\medskip For $a\in \mathbb{N} \cup\{\infty\}$ we will define a
functor:

\begin{equation}\label{VZ0e-1c}
  \mathfrak{S}_a\text{-windows} \longrightarrow \;
  \text{Dieudonn\'e displays over} \;  R/p^aR.
\end{equation}
Let $(Q,\phi)\equiv (P,Q,F)$ be a $\mathfrak{S}_a$-window. Let
$F_1:Q\to P$ be as in section 2. To $(P,Q,F)$ we will associate a
Dieudonn\'e display $(P',Q',F',F_1')$ over $R/p^aR$. Let $P' :=
\hat{W}(R/p^aR) \otimes_{\kappa_a,\mathfrak{S}_a} P$. Let $Q'$ be the
kernel of the natural $\hat{W}(R/p^aR)$-linear epimorphism:
\begin{displaymath}
  P'= \hat{W}(R/p^aR) \otimes_{\kappa_a,\mathfrak{S}_a} P \twoheadrightarrow P/Q.
\end{displaymath}
We define $F' : P' \rightarrow P'$ as the canonical $F$-linear
extension of $F$. We define $F'_1 : Q' \rightarrow P'$ by the rules:
\begin{displaymath}
  \begin{array}{lcrr}
    F'_1(\xi \otimes y) & = & ~^F\xi \otimes \tau F_1(y),& \mbox{for} \; \xi \in
    \hat{W}(R/p^aR), \; y \in Q\\
    F'_1(~^V\xi \otimes x) & = & \xi \otimes F(x), & \mbox{for} \; \xi \in
    \hat{W}(R/p^aR), \; x \in P.\\
  \end{array}
\end{displaymath} 
Using a normal decomposition of $(P,Q,F)$, one checks that
$(P',Q',F',F_1')$ is a Dieudonn\'e display over $R/p^aR$.

Since the category of Dieudonn\'e displays over $R/p^aR$ is equivalent
to the category of $p$-divisible groups over $R/p^aR$ (see \cite{ZDD})
we obtain from (\ref{VZ0e-1c}) a functor
\begin{equation}\label{VZ0-1d}
  \mathfrak{S}_a\text{-windows} \longrightarrow \;
  p\text{-divisible groups over}\;  R/p^aR.  
\end{equation}
In particular, for the $p$-divisible group $G$ associated to
$(Q,\phi)\equiv (P,Q,F)$ we have identifications of $R/p^aR$-modules
\begin{equation}\label{VZ0-1l}
  \Lie(G)=P'/Q'=P/Q=\text{Coker}(\phi).
\end{equation}
 
\smallskip
 
In $\mathfrak{S}_1$ the elements $E$ and $p$ differ by a
unit. Therefore the notion of a $\mathfrak{S}_1$-window is the same as
that of a Dieudonn\'e $\mathfrak{S}_1$-window over $R/pR$ introduced
in \cite{ZW}, Definition 2. By Theorem 6 (or 3.2) of loc. cit. we get
that the functor (\ref{VZ0-1d}) is an equivalence of categories in the
case $a=1$. We would like to mention that the contravariant analogue
of this equivalence for $a=1$ also follows from \cite{dJ}, Theorem of
Introduction and Proposition 7.1.

The faithfulness of the functors (\ref{VZ0e-1c}) and (\ref{VZ0-1d}) follows from the mentioned equivalence in the case $a=1$ and from the following rigidity
property:

\begin{lemma}
  Let $a\geq 1$ be a natural number.  Let $\mathcal{P} = (P,Q,F)$ and
  $ \mathcal{P}' = (P',Q',F')$ be $\mathfrak{S}_{ap}$-windows. By base
  change we obtain windows $\bar{\mathcal{P}}$ and
  $\bar{\mathcal{P}}'$ over $\mathfrak{S}_{a}$. Then the natural map

\begin{displaymath}
  Hom_{\mathfrak{S}_{ap}} (\mathcal{P}, \mathcal{P}') \rightarrow
  Hom_{\mathfrak{S}_{a}} (\bar{\mathcal{P}}, \bar{\mathcal{P}}') 
\end{displaymath}
is injective.
\end{lemma} {\bf Proof:}
Let $\alpha : \mathcal{P} \rightarrow \mathcal{P}'$ be a morphism, which
induces $0$ over $\mathfrak{S}_a$. We have $\alpha(P) \subset
u^{ae}P'$. To prove that $\alpha = 0$ it is enough to show that $\alpha
(F_1y) = 0$ for each $y \in Q$. We have $\alpha(y) \in u^{ae}P' \cap
Q'$. We choose a normal decomposition $P' = J' \oplus L'$ and we write
$\alpha(y) = j' + l'$. Then we have $j' \in u^{ae}J' \cap EJ' =
u^{ae}EJ'$ and $l' \in u^{ae}L'$. In $\mathfrak{S}_{ap}$ we have
$\sigma (u^{ae}) = u^{ape} = 0$. We conclude that $F'_1j' = 0$ and
$F'_1l' = 0$. Finally we obtain $\alpha (F_1y) = F'_1(\alpha y) =
0$. \endproof

\begin{lemma}\label{VZ2l}
The functors (\ref{VZ0e-1c}) and
  (\ref{VZ0-1d}) are essentially surjective on objects.
\end{lemma}
{\bf Proof}: We will first prove the lemma for $a\in \mathbb{N}$. We will
use induction on $a\in \mathbb{N}$. We 
already know that this is true for $a = 1$. The inductive passage
from $a$ to $a+1$ goes as follows. It suffices to consider the case of
the functor (\ref{VZ0e-1c}).
 
\smallskip
 
Let $\tilde{\mathcal{P}}' = (\tilde{P}', \tilde{Q}', \tilde{F}',
\tilde{F}_1')$ be a Dieudonn\'e display over $R/p^{a+1}R$. We denote
by $\mathcal{P}'=(P',Q',F',F'_1)$ its reduction over $R/p^aR$. Then we
find by induction a $\mathfrak{S}_a$-window $\mathcal{P}$ which is
mapped to $\mathcal{P}'$ by the functor (\ref{VZ0e-1c}). We lift
$\mathcal{P}$ to a $\mathfrak{S}_{a+1}$-window
$\tilde{\mathcal{P}}=(\tilde{Q},\tilde\phi)\equiv
(\tilde{P},\tilde{Q},\tilde{F})$, cf. end of section 2. Let $\tilde
F_1:\tilde Q\to\tilde P$ be obtained from $\tilde F$ as in section
2.

We apply to $\tilde{\mathcal{P}}$ the functor (\ref{VZ0e-1c}) and we
obtain a Dieudonn\'e display $\tilde{\mathcal{P}}''=(\tilde{P}'',
\tilde{Q}'', \tilde{F}'',\tilde{F}_1'')$ over $R/p^{a+1}R$. By
\cite{ZDD}, Theorem 3 we can identify
\begin{equation}\label{VZ0e-1e}
  (\tilde{P}', \tilde{F}', \Phi_1)=(\hat{W}(R/p^{a+1}R)
  \otimes_{\kappa_{a+1},\mathfrak{S}_{a+1}} \tilde{P} = \tilde{P}'',
  \tilde{F}'', \Phi_1). 
\end{equation}
Here $\Phi_1: \breve{Q}' \rightarrow \tilde{P}'$ is a Frobenius linear
map from the inverse image $\breve{Q}'$ of $Q'$ in
$\tilde{P}'=\tilde{P}''$ which extends both $\tilde{F}_1'$ and
$\tilde{F}_1''$ and which satisfies the identity
$\Phi_1([p^a]\tilde{P}') = 0$ (this identity is due to the fact that
we use the trivial divided power structure on the kernel of the
epimorphism $R/p^{a+1}R \twoheadrightarrow R/p^aR$).

The composite map:
\begin{displaymath}
  \tilde{Q} \overset{\tilde{F}_1}{\longrightarrow } \tilde{P} \rightarrow
  \tilde{P}' \overset{\tau}{\rightarrow} \tilde{P}', 
\end{displaymath}  
coincides with the composite map
\begin{displaymath}
  \tilde{Q} \rightarrow \breve{Q}' \overset{\Phi_1}{\longrightarrow }
  \tilde{P}'. 
\end{displaymath}

We define $\tilde{Q}^{*} \subset \tilde{P}$ as the inverse image of
the natural map $\tilde{P} \rightarrow \tilde{P}'/\tilde{Q}'$ deduced
from the identity $\tilde{P}'=\hat{W}(R/p^{a+1}R)
\otimes_{\kappa_{a+1},\mathfrak{S}_{a+1}} \tilde{P}$.  The images of
$\tilde{Q}$ and $\tilde{Q}^{*}$ by the canonical map $\tilde{P}
\rightarrow P$ are the same.  Therefore for each $y^{*} \in
\tilde{Q}^{*}$ there exists an $y \in \tilde{Q}$ such that we have
$y^{*} = y + u^{ae}x$ for some $x \in \tilde{P}$. Since $\tilde{F}
(u^{ae}x) = 0$ we conclude that $\tilde F(y^{*}) = \tilde F(y) \in
\sigma(E)\cdot\tilde{P}$. This proves that $\tilde{\mathcal{P}}^{*} =
(\tilde{P}, \tilde{Q}^{*}, \tilde{F})$ is a
$\mathfrak{S}_{a+1}$-window which lifts the $\mathfrak{S}_a$-window
$\mathcal{P}$. Let $\tilde F_1^{*}:\tilde Q^{*}\to\tilde P$ be
obtained from $\tilde F$ as in section 2.

We claim that the image of $\tilde{\mathcal{P}}^{*}$ via the functor
(\ref{VZ0e-1c}) coincides with the Dieudonn\'e display
$\tilde{\mathcal{P}}'$. For this we have to show that the composite
map

\begin{displaymath}
  \tilde{Q}^{*} \overset{\tilde{F}_1^{*}}{\longrightarrow } \tilde{P} \rightarrow
  \tilde{P}' \overset{\tau}{\rightarrow} \tilde{P}', 
\end{displaymath}  
coincides with the composite map
\begin{displaymath}
  \tilde{Q}^{*} \rightarrow \breve{Q}' \overset{\Phi_1}{\longrightarrow }
  \tilde{P}'. 
\end{displaymath}
This follows again from the decomposition $y^{*} = y + u^{ae}x$ and
the facts that: (i) we have $\tilde F_1^*(y^{*})=\tilde F_1(y)$ (as we
have $\tilde F(y^{*})=\tilde F(y)$) and (ii) the image of $u^{ae}x$ in
$\breve{Q}'$ is mapped to zero by $\Phi_1$. We conclude that
$\tilde{\mathcal{P}}'$ is in the essential image of the functor
(\ref{VZ0e-1c}). This ends the induction.

The fact that the lemma holds even if $a=\infty$ follows from the above induction via a natural limit process. \endproof

\bigskip

\section{Extending morphisms between $\mathfrak{S}_1$-windows}

In this section we prove an extension result for an isomorphism
between $\mathfrak{S}_1$-windows. We begin by considering for $a\in
\mathbb{N} \cup\{\infty\}$ an extra $W(k)$-algebra:
\begin{displaymath}
  \mathcal{T}_a :=  \mathfrak{S}[[v]]/(pv - u^e, v^a)
\end{displaymath}
with the convention that $v^{\infty}:=0$. This ring is without
$p$-torsion. It is elementary to check that the canonical ring
homomorphism 
\begin{displaymath}
   \mathfrak{S}_a \rightarrow \mathcal{T}_a
\end{displaymath}
is injective. For $a=1$ this is an isomorphism $\mathfrak{S}_1 \cong
\mathcal{T}_1 $. 

We set $\mathcal{T} =
\mathcal{T}_{\infty}$.  In $\mathcal{T}_a$ we have
$E = p(v + \epsilon)$ and thus the elements $p$ and $E$ differ by a
unit. We have an isomorphism:
\begin{displaymath}
  \mathcal{T}_a/p\mathcal{T}_a \cong (R/pR)[[v]]/(v^a).
\end{displaymath} 
We extend the Frobenius endomorphism $\sigma $ to $\mathcal{T}_a$ by
the rule:
\begin{displaymath}
  \sigma (v) = u^{e(p-1)}v = p^{p-1}v^p
\end{displaymath}

We note that the endomorphism $\sigma $ on $\mathcal{T}_a$ no longer
induces the Frobenius modulo $p$. But the notion of a window over
$\mathcal{T}_a$ still makes sense as follows.

\begin{df} A window over $\mathcal{T}_a$ is a triple $(P,Q,F)$,
  where $P$ is a free $\mathcal{T}_a$-module, $Q$ is a
  $\mathcal{T}_a$-submodule of $P$ such that $P/Q$ is a free
  $\mathcal{T}_a/p\mathcal{T}_a$-module, and $F: P \rightarrow P$ is a
  $\sigma$-linear endomorphism. We require that $F(Q) \subset pP$ and
  that this subset generates $pP$ as a $\mathcal{T}_a$-module.  
\end{df}
We define a $\sigma$-linear map $F_1: Q \rightarrow P$ by
$pF_1(y) = F(y)$ for $y\in Q$. Its linearisation $F_1^{\sharp}$ is an
isomorphism. Taking the composite of the 
inclusion $Q \subset P$ with $(F_1^{\sharp})^{-1}$ we obtain a
$\mathcal{T}_a$-linear map
\begin{displaymath}
  \phi : Q \rightarrow Q^{(\sigma)},
\end{displaymath}
whose cokernel is a free $\mathcal{T}_a/p\mathcal{T}_a$-module.

If we start with a triple $(P,Q,F)$ as in Lemma \ref{Tripel} an
tensor it with $\mathcal{T}_a \otimes_{\mathfrak{S}_a} $ we obtain a
window over $\mathcal{T}_a$. 

A window over $\mathcal{T}_a$ is not a window in sense of \cite{ZW} because
$\sigma$ on $\mathcal{T}_a/p\mathcal{T}_a$ is not the Frobenius
endomorphism. We have still the following lifting property.

\begin{prop}\label{VZ3l}
  Let $(Q_1,\phi_1)$ and $(Q_2,\phi_2)$ be two Breuil windows relative
  to $\mathfrak{S} \rightarrow R$. Let $(\breve{Q}_1, \breve{\phi}_1)$
  and $(\breve{Q}_2, \breve{\phi}_2)$ be the $\mathfrak{S}_1$-windows
  which are the reduction modulo $u^e$ of $(Q_1,\phi_1)$ and
  $(Q_2,\phi_2)$ (respectively). Let $\breve{\alpha}: \breve{Q}_1
  \rightarrow \breve{Q}_2$ be an isomorphism of windows relative to
  $\mathfrak{S}_1 \rightarrow R/pR$ i.e., a $\mathfrak{S}_1$-linear
  isomorphism such that we have $\breve{\phi}_2\circ
  \breve{\alpha}=(1\otimes \breve{\alpha})\circ \breve{\phi}_1$.

  Then there exists a unique isomorphism
  \begin{displaymath}
    \alpha: \mathcal{T} \otimes_{\mathfrak{S}} Q_1 \rightarrow
    \mathcal{T} \otimes_{\mathfrak{S}} Q_2 
  \end{displaymath}
  which commutes in the natural sense with $\phi_1$ and $\phi_2$ and
  which lifts $\breve{\alpha} $ with respect to the
  $\mathfrak{S}$-epimorphism $\mathcal{T} \rightarrow \mathfrak{S}_1$
  that maps $v$ to $0$.
\end{prop} {\bf Proof:} We choose a normal decomposition $\breve{Q}_1
= \breve{L}_1 \oplus \breve{J}_1$. Applying $\breve{\alpha}$ we obtain
a normal decomposition $\breve{Q}_2 = \breve{L}_2 \oplus \breve{J}_2$.
We lift these normal decompositions to $\mathfrak{S}$:
\begin{displaymath}
  Q_1 = J_1 \oplus L_1\;\; \qquad \mbox{and}\;\; \qquad Q_2 = J_2 \oplus L_2.
\end{displaymath}
We find an isomorphism $\gamma : Q_1 \rightarrow Q_2$ which lifts
$\breve{\alpha}$ and such that $\gamma (L_1) = L_2$ and $\gamma (J_1)
= J_2$. We identify the modules $Q_1$ and $Q_2$ via $\gamma$ and we
write:
\begin{displaymath}
  Q = Q_1 = Q_2, \quad J = J_1 = J_2, \quad L = L_1 = L_2.
\end{displaymath}

We choose a $\mathfrak{S}$-basis $\{e_1, \ldots, e_d\}$ for $J$ and a
$\mathfrak{S}$-basis $\{e_{d+1}, \ldots e_r\}$ for $L$. Then
$\{1\otimes e_1, \ldots, 1 \otimes e_r\}$ is a $\mathfrak{S}$-basis
for $Q^{(\sigma)}$. For $i\in\{1,2\}$ we write $\phi_i: Q \rightarrow
Q^{(\sigma)}$ as a matrix with respect to the mentioned
$\mathfrak{S}$-bases. It follows from the properties of a normal
decomposition that this matrix has the form:
\begin{displaymath}
  A_i 
  \left(
    \begin{array}{rr}
      E\cdot I_d & 0\\
      0    & I_c
    \end{array}
  \right),
\end{displaymath} 
where $A_1$ and $A_2$ are invertible matrices in $\GL_r(\mathfrak{S})$
and where $c:=r-d$. By the construction of $\gamma$, the
$\mathfrak{S}$-linear maps $\phi_1$ and $\phi_2$ coincide modulo
$(u^e)$. From this and the fact that $E$ modulo $(u^e)$ is a non-zero
divisor of $\mathfrak{S}/(u^e)$, we get that we can write
\begin{equation}\label{VZ4e}
  (A_2)^{-1} A_1 = I_r + u^eZ,
\end{equation}
where $Z \in M_r(\mathfrak{S})$.  We set
\begin{displaymath}
  C := \left(
    \begin{array}{rr}
      E\cdot I_d & 0\\
      0    & I_c
    \end{array}
  \right).
\end{displaymath}

To find the isomorphism $\alpha$ is the same as to find a matrix $X
\in \GL_r(\mathcal{T})$ which solves the equation

\begin{equation}\label{VZ5e}
  A_2 C X = \sigma(X) A_1 C 
\end{equation}
and whose reduction modulo the ideal $(v)$ of $\mathcal{T}$ is the
matrix representation of $\breve{\alpha}$ i.e., it is the identity
matrix. Therefore we set:
\begin{equation}\label{VZ6e}
  X = I_r + vY
\end{equation}
for a matrix $Y \in M_r(\mathcal{T})$.

As $E/p=v+\epsilon$ is a unit of $\mathcal{T}$, the matrix $pC^{-1}$
has coefficients in $\mathcal{T}$.

From the equations (\ref{VZ5e}) and (\ref{VZ6e}) we obtain the
equation:
\begin{displaymath}
  p(I_r + vY) = pC^{-1}A_2^{-1}(I_r + \sigma(v)\sigma(Y))A_1C.
\end{displaymath}
We insert (\ref{VZ4e}) in this equation. With the notation $D :=
pC^{-1}ZC$ we find:
\begin{displaymath}
  u^eY = u^eD + (u^{ep}/p) pC^{-1}A_2^{-1}\sigma(Y)A_1C.
\end{displaymath}
Since $u^e$ is not a zero divisor in $\mathcal{T}$ we can write:
\begin{equation}\label{VZ7e}
  Y - (u^{e(p-1)}/p) pC^{-1}A_2^{-1}\sigma(Y)A_1C = D.
\end{equation}
We have $(u^{e(p-1)}/p)=vu^{e(p-2)}$. Thus the $\sigma$-linear
operator
$$\Psi (Y) = (u^{e(p-1)}/p)
pC^{-1}A_2^{-1}\sigma(Y)A_1C$$ on the $\mathcal{T}$-module
$M_r(\mathcal{T})$ is topologically nilpotent. Therefore the equation
(\ref{VZ7e}) has a unique solution $Y=\sum_{n=0}^{\infty} \Psi(D)\in
M_r(\mathcal{T})$. Therefore $X = I_r + vY\in \GL_r(\mathcal{T})$
exists and is uniquely determined. \endproof

\section{Proof of Theorem \ref{VZ1t}}

In this section we prove Theorem \ref{VZ1t}. In section 3 we
constructed a functor (cf. (\ref{VZ0-1d}) with $a=\infty$)
\begin{equation}\label{VZ0-1s}
  \mbox{Breuil} \; \mbox{windows} \; \mbox{relative} \; \mbox{to}\; \mathfrak{S}\to R\;\;\;\longrightarrow \;\;\;
  p\mbox{-divisible} \; \mbox{groups} \; \mbox{over} \;  R
\end{equation}
which is faithful and which (cf. Lemma \ref{VZ2l}) is essentially surjective on objects.
Thus to end the proof of Theorem \ref{VZ1t}, it suffices to show that
this functor is essentially surjective on isomorphisms (equivalently,
on morphisms). This will be proved in Lemma \ref{VZ4l} below. For the
proof of Lemma \ref{VZ4l}, until the end we take $a\in \mathbb{N}$ and
we will begin by first listing some basic properties of the rings
$\mathfrak{S}_a$ and $\mathfrak{T}_a$.

\subsection{On $\mathfrak{S}_a$ and $\mathfrak{T}_a$} 
There exists a canonical homomorphism
\begin{displaymath}
  \mathfrak{S}_a \rightarrow R/p^aR
\end{displaymath} 
whose kernel is the principal ideal of $\mathfrak{S}_a$ generated by
$E$ modulo $(u^{ae})$. Let
\begin{displaymath}
  \mathcal{S}_a \subset \mathfrak{S}_a \otimes \mathbb{Q}
\end{displaymath}
be the subring generated by all elements $u^{ne}/n!$ over
$\mathfrak{S}_a$ with $n\in\{0,\ldots,a\}$. Then $\mathcal{S}_a$ is a
$p$-adic ring without $p$-torsion. There exists a commutative diagram
\begin{displaymath}
  \xymatrix{
    \mathfrak{S}_a \ar[rd] \ar[d] &\\
    \mathcal{S}_a \ar[r] & R/p^{\nu(a)}R,
  }
\end{displaymath}
where $\nu(a) := \inf \{\ord_p(\frac{p^n}{n!})|\; \mbox{for} \; n \geq
a\}.$ With the notation of (\ref{VZ0e}), the horizontal homomorphism
maps $u^{ne}/n!$ to $(p^n/n!) (-\epsilon)^n$. By \cite{ZW}, Theorem 6
or 3.2, the epimorphism $\mathcal{S}_a \rightarrow R/p^{\nu(a)}R$ is a
frame which classifies $p$-divisible groups over $R/p^{\nu(a)}R$. In
what follows, by a $\mathcal{S}_a$-window we mean a Dieudonn\'e
$\mathcal{S}_a$-window over $R/p^{\nu(a)}R$ in the sense of \cite{ZW},
Definition 2.

Both $\mathcal{S}_a$ and $\mathcal{T}_a$ are subrings of
$\mathfrak{S}_a \otimes \mathbb{Q}$. Because of the identity
\begin{displaymath}
  \frac{u^{ne}}{n!} = \frac{p^n}{n!} \left(\frac{u^e}{p}\right)^n =
  \frac{p^n}{n!} v^n
\end{displaymath}
we have an inclusion of $W(k)$-algebras:
\begin{displaymath}
  \mathcal{S}_a \subset \mathcal{T}_a.
\end{displaymath}
The Frobenius endomorphism $\sigma$ of $\mathfrak{S}$ induces a
homomorphism
\begin{displaymath}
  \mathfrak{S}_a \rightarrow \mathfrak{S}_{pa}.
\end{displaymath}
It maps the subalgebra $\mathcal{T}_a \subset \mathfrak{S}_a \otimes
\mathbb{Q}$ to $\mathcal{S}_{pa} \subset \mathfrak{S}_{pa} \otimes
\mathbb{Q}$. We denote this homomorphism by
\begin{displaymath}
  \tau_a : \mathcal{T}_a \rightarrow  \mathcal{S}_{pa}.
\end{displaymath} 

\begin{lemma}\label{VZ4l}
  Let $(Q_1,\phi_1)$ and $(Q_2,\phi_2)$ be two Breuil windows relative
  to $\mathfrak{S} \rightarrow R$. Let $G_1$ and $G_2$ be the
  corresponding $p$-divisible groups over $R$, cf. the functor
  (\ref{VZ0-1s}). Let $\gamma :G_1 \rightarrow G_2$ be an
  isomorphism. Then there exists a unique morphism (automatically
  isomorphism) $\alpha_0:(Q_1,\phi_1)\to (Q_2,\phi_2)$ of Breuil
  windows relative to $\mathfrak{S} \rightarrow R$ which maps to
  $\gamma$ via the functor (\ref{VZ0-1s}).
\end{lemma} {\bf Proof:} For $i\in\{1,2\}$ let
$(\breve{Q}_i,\breve{\phi}_i)\equiv
(\breve{P}_i,\breve{Q}_i,\breve{F}_i)$ be the $\mathfrak{S}_1$-window
which is the reduction of $(Q_i,\phi_i)\equiv (P_i,Q_i,F_i)$ modulo
$u^e$. Since the category of $\mathfrak{S}_1$-windows is equivalent to the category of $p$-divisible
groups over $R/pR$, the reduction of $\gamma$ modulo $p$ is induced
via the functor (\ref{VZ0-1d}) by an isomorphism
$\breve{\alpha}:(\breve{Q}_1,\breve{\phi}_1)\to
(\breve{Q}_2,\breve{\phi}_2)$. Due to Proposition \ref{VZ3l}, the last
isomorphism extends to an isomorphism
$\alpha:\mathcal{T}\otimes_{\mathfrak{S}} (Q_1,\phi_1)\to
\mathcal{T}\otimes_{\mathfrak{S}} (Q_2,\phi_2)$ of windows over
$\mathcal{T}$. Here and in what follows we identify $\mathcal{T}_a\otimes_{\mathfrak{S}} (Q_i,\phi_i)\equiv \mathcal{T}_a\otimes_{\mathfrak{S}} (P_i,Q_i,F_i)$ and therefore we  will refer to $\mathcal{T}_a\otimes_{\mathfrak{S}} (Q_i,\phi_i)$ as a window over $\mathcal{T}_a$ (here $a$ can be also $\infty$). As in the proof of Proposition \ref{VZ3l} we can
identify normal decompositions $Q_1=J_1\oplus L_1=J_2\oplus L_2=Q_2$
and we can represent the mentioned isomorphism of windows over
$\mathcal{T}$ by an invertible matrix $X \in \GL_r(\mathcal{T})$.  Let
$X_a\in \GL_r(\mathcal{T}_a)$ be the reduction of $X$ modulo $v^a$.

The matrix $X$ has the following crystalline interpretation.  The
epimorphism $\mathcal{T}_a \twoheadrightarrow (R/pR)[[v]]/(v^a)$ is a
pd-thickening. (We emphasize that it is not a frame in the sense of
\cite{ZW}, Definition 1 because $\sigma$ modulo $p$ is not the
Frobenius endomorphism of $\mathcal{T}_a/p\mathcal{T}_a $.)  We have a
morphism of pd-thickenings
\begin{equation}\label{VZe8}
  \xymatrix{
    \mathcal{S}_a \ar[r] \ar[d] & \mathcal{T}_a \ar[d] \\
    R/p^{\nu(a)}R \ar[r] & (R/pR)[[v]]/(v^a).
  }
\end{equation}
By the crystal associated to a $p$-divisible group $\square$ over
$R/p^{\nu(a)}R$ we mean the Lie algebra crystal of the universal
vector extension crystal of $\square$ as defined in [M]. The crystal
of $G_i$ evaluated at the pd-thickening $\mathcal{S}_a \rightarrow
R/p^{\nu(a)}R$ coincides in a functorial way with $\mathcal{S}_a
\otimes_{\mathfrak{S}_a} Q_i^{(\sigma)}=\mathcal{S}_a
\otimes_{\mathfrak{S}_a} P_i$ (if $G_i$ is a formal $p$-divisible
group, this follows from either \cite{ZDFG}, Theorem 6 or from
\cite{ZW}, Theorem 1.6; the general case follows from \cite{L}).  Let
$\breve{G}_{i,a}$ be the push forward of $G_i$ via the canonical
homomorphism $R \rightarrow (R/pR)[[v]]/(v^a)$. The diagram
(\ref{VZe8}) shows that $\mathcal{T}_a \otimes_{\mathfrak{S}}
Q_i^{(\sigma)}=\mathcal{T}_a \otimes_{\mathfrak{S}} P_i$ is the
crystal of $\breve{G}_{i,a}$ evaluated at the pd-thickening
$\mathcal{T}_a \rightarrow (R/pR)[[v]]/(v^a)$.  The isomorphism
$\gamma: G_1 \rightarrow G_2$ induces an isomorphism $\alpha_a$ of
$\mathcal{S}_a$-windows which is defined by a $\mathcal{S}_a$-linear
isomorphism $\mathcal{S}_a\otimes_{\mathfrak{S}} P_1\to
\mathcal{S}_a\otimes_{\mathfrak{S}} P_2$. Via base change of
$\alpha_a$ through the diagram (\ref{VZe8}), we get an isomorphism
$\beta_a:\mathcal{T}_a\otimes_{\mathfrak{S}} (Q_1,\phi_1)\to
\mathcal{T}_a\otimes_{\mathfrak{S}} (Q_2,\phi_2)$ of windows over
$\mathcal{T}_a$ defined by an isomorphism
$\mathcal{T}_a\otimes_{\mathfrak{S}} (P_1,Q_1,F_1)\to
\mathcal{T}_a\otimes_{\mathfrak{S}} (P_2,Q_2,F_2)$. We note that after
choosing a normal decomposition a window is simply an invertible
matrix (to be compared with end of section 2) and base change applies
to the coefficients of this matrix the homomorphism $\mathcal{S}_a
\rightarrow \mathcal{T}_a $. The system of isomorphisms $\beta_a$
induces in the limit an isomorphism of windows over $\mathcal{T}$:
\begin{equation}
  \beta: \mathcal{T} \otimes_{\mathfrak{S}} (Q_1,\phi_1) \rightarrow
  \mathcal{T} \otimes_{\mathfrak{S}} (Q_2,\phi_2) .
\end{equation}  
We continue the base change (\ref{VZe8}) using the following diagram:
\begin{equation}
  \xymatrix{
    \mathcal{T}_a \ar[r] \ar[d] & \mathfrak{S}_1 \ar[d] \\
    (R/pR)[[v]]/(v^a) \ar[r] & R/pR.
  }
\end{equation}
The vertical arrows are thickenings with divided powers.  From
$\alpha_a$ we obtain by base change the isomorphism $\breve{\alpha}$
since windows associated to $p$-divisible groups commute with base
change. But this shows that the isomorphism $\beta$ coincides with the
isomorphism $\alpha: \mathcal{T} \otimes_{\mathfrak{S}} (Q_1,\phi_1)
\rightarrow \mathcal{T} \otimes_{\mathfrak{S}} (Q_2,\phi_2)$, cf. the
uniqueness part of Proposition \ref{VZ3l}.

We will show that the assumption that $X$ has coefficients in
$\mathfrak{S}$ implies the Lemma. This assumption implies that
$\alpha$ is the tensorization with $\mathcal{T}$ of an isomorphisms
$\alpha_0:(Q_1,\phi_1)\to (Q_2,\phi_2)$. Let $\gamma_0:G_1\to G_2$ be
the isomorphism associated to the isomorphism $\alpha_0$ via the
functor (\ref{VZ0-1s}). As the functor (\ref{VZ0-1d}) is an
equivalence of categories for $a=1$, we get that $\gamma_0$ and
$\gamma$ coincide modulo $p$. Therefore $\gamma=\gamma_0$ and thus the
Lemma holds.

Thus to end the proof of the Lemma it suffices to show by induction on
$a\in \mathbb{N}$ that the matrix $X_a$ has coefficients in
$\mathfrak{S}_a$. The case $a =1$ is clear. The inductive passage from
$a$ to $a+1$ goes as follows. We can assume that $X_a$ has
coefficients in $\mathfrak{S}_a$. Therefore the invertible matrix
$\tau_a(X_a)\in \GL_r(\tau(\mathcal{T}_a))\subset
\GL_r(\mathcal{S}_{pa})$ defines a $\mathcal{S}_{pa}$-linear
isomorphism

\begin{displaymath}
  \mathcal{S}_{pa} \otimes_{\mathfrak{S}} P_1=\mathcal{S}_{pa}
  \otimes_{\mathfrak{S}} Q_1^{(\sigma)} \rightarrow 
  \mathcal{S}_{pa} \otimes_{\mathfrak{S}} Q_2^{(\sigma)} =
  \mathcal{S}_{pa} \otimes_{\mathfrak{S}} P_2 
\end{displaymath}
which respects the Hodge filtration i.e., it is compatible with the
$R/p^{\nu(pa)}R$-linear map $\Lie G_{1,R/p^{\nu(pa)}R} \rightarrow
\Lie G_{2,R/p^{\nu(pa)}R}$ induced by $\gamma$.

Since $\tau_a(X_a)$ has coefficients in $\mathfrak{S}_{pa}$, we obtain
a commutative diagram of $\mathfrak{S}_{pa}$-modules:
\begin{displaymath}
  \xymatrix{
    \mathfrak{S}_{pa} \otimes_{\mathfrak{S}} Q_1^{(\sigma)} \ar[r] \ar[d] 
    & \Lie G_{1,R/p^{\nu(pa)}R} \ar[d] \\
    \mathfrak{S}_{pa} \otimes_{\mathfrak{S}} Q_2^{(\sigma)} \ar[r]  
    & \Lie G_{1,R/p^{\nu(pa)}R}.
  }
\end{displaymath}  
Since $\nu(pa) \geq a+1$, from (\ref{VZ0-1l}) we obtain a commutative
diagram
\begin{equation}\label{VZ9e}
  \xymatrix{
    0 \ar[r] & \mathfrak{S}_{a+1} \otimes_{\mathfrak{S}} Q_1 \ar[r] &
    \mathfrak{S}_{a+1} \otimes_{\mathfrak{S}} Q_1^{(\sigma)} \ar[r] \ar[d]  
    & \Lie G_{1,R/p^{(a+1)}R} \ar[d] \\
    0 \ar[r] & \mathfrak{S}_{a+1} \otimes_{\mathfrak{S}} Q_2 \ar[r] &
    \mathfrak{S}_{a+1} \otimes_{\mathfrak{S}} Q_2^{(\sigma)} \ar[r]   
    & \Lie G_{1,R/p^{(a+1)}R} 
  }  
\end{equation}
with exact rows. The left vertical arrow is induced by
$\sigma(X_{a+1})$.  On the kernels of the horizontal maps we obtain a
$\mathfrak{S}_{a+1}$-linear isomorphism
\begin{equation}\label{VZ10e}
  \mathfrak{S}_{a+1} \otimes_{\mathfrak{S}} Q_1 \rightarrow
  \mathfrak{S}_{a+1} \otimes_{\mathfrak{S}} Q_2.
\end{equation}
As $E$ is not a zero divisor in $\mathcal{T}_{a+1}$, by tensoring the
short exact sequences of (\ref{VZ9e}) with $\mathcal{T}_{a+1}$ we get
short exact sequences of $\mathcal{T}_{a+1}$-modules. This implies
that the tensorization of the $\mathfrak{S}_{a+1}$-linear isomorphism
(\ref{VZ10e}) with $\mathcal{T}_{a+1}$ is given by the matrix
$X_{a+1}$. Therefore $X_{a+1}$ has coefficients in
$\mathfrak{S}_{a+1}$. This completes the induction and thus the proofs
of the Lemma and of Theorem \ref{VZ1t}.\endproof

\section{Breuil modules}

To prove Theorem \ref{VZ2t} we first need the following basic result
on Breuil modules relative to $\mathfrak{S}\rightarrow R$.
 
\begin{prop}\label{P2}
  Let $(M, \varphi)$ be a Breuil module relative to $\mathfrak{S}
  \rightarrow R$. Then the following four properties hold:

  \begin{enumerate}
  \item[(i)] The $\mathfrak{S}$-linear map $\varphi$ is injective.
  \item[(ii)] There exists a short exact sequence $0\to (Q',\phi')\to
    (Q,\phi)\to (M,\varphi)\to 0$, where $(Q',\phi')$ and $(Q,\phi)$
    are Breuil windows relative to $\mathfrak{S} \rightarrow R$.
  \item[(iii)] If $(M,\varphi)\to (\tilde{M},\tilde{\varphi})$ is a
    morphism of Breuil modules relative to $\mathfrak{S}\rightarrow
    R$, then it is the cokernel of a morphism between two exact
    complexes $0\to (Q',\phi')\to (Q,\phi)$ and $0\to
    (\tilde{Q}',\tilde{\phi}')\to (\tilde{Q},\tilde{\phi})$ of Breuil
    windows relative to $\mathfrak{S} \rightarrow R$.
  \item[(iv)] The quotient $M^{(\sigma)} /\varphi (M)$ is an
    $R$-module of projective dimension at most one.
  \end{enumerate}
\end{prop} {\bf Proof:}
Let $(p):=p\mathfrak{S}$; it is a principal prime ideal of
$\mathfrak{S}$. Then $\sigma$ induces an endomorphism of the discrete
valuation ring $\mathfrak{S}_{(p)}$ which fixes $p$. Thus the length
of a $\mathfrak{S}_{(p)}$-module remains unchanged if tensored by
$\sigma_{(p)}:\mathfrak{S}_{(p)} \rightarrow \mathfrak{S}_{(p)}$. One
easily checks that
\begin{displaymath}
  (M_{(p)})^{(\sigma)} = \mathfrak{S}_{(p)} \otimes_{\sigma_{(p)},\mathfrak{S}_{(p)}}
  M_{(p)} \cong (M^{(\sigma)})_{(p)}.
\end{displaymath}

Let $x, p$ be a regular sequence in $\mathfrak{S}$. As the
$\mathfrak{S}$-module $M$ has projective dimension at most one and as
$M$ is annihilated by a power of $p$, the multiplication by $x$ is a
$\mathfrak{S}$-linear monomorphism $x : M \hookrightarrow M$. Since no
element of $\mathfrak{S}\setminus p\mathfrak{S}$ is a zero divisor in
$M$, we conclude that $M \subset M_{(p)}$. The $\mathfrak{S}$-linear
map $\varphi: M \rightarrow M^{(\sigma)}$ becomes an epimorphism when
tensored with $\mathfrak{S}_{(p)}$. We obtain an epimorphism of
$\mathfrak{S}_{(p)}$-modules of the same length:
\begin{displaymath}
  \varphi_{(p)} : M_{(p)} \rightarrow (M_{(p)})^{(\sigma)}.
\end{displaymath}
As $M$ is a finitely generated $\mathfrak{S}$-module annihilated by a
power of $p$, the $\mathfrak{S}_{(p)}$-module $M_{(p)}$ has finite
length. From the last two sentences we get that $\varphi_{(p)}$ is
injective. Thus $\varphi$ is also injective i.e., (i) holds.

We consider free $\mathfrak{S}$-modules $J$ and $L$ of finite ranks
and a $\mathfrak{S}$-linear epimorphism
\begin{displaymath}
  J \oplus L \overset{\tau}{\longrightarrow } M^{(\sigma)}
\end{displaymath}  
which maps the free $\mathfrak{S}$-submodule $EJ \oplus L$
surjectively to $\varphi (M)$. Let $\tau_1: J \oplus L \rightarrow M$
be the unique $\mathfrak{S}$-linear map such that we have a
commutative diagram
\begin{displaymath}
  \begin{CD}
    EJ \oplus L @>>> \Image (\varphi)\\
    @A{E\id_J + \id_L}AA  @AA {\varphi}A\\
    J \oplus L @>>> M
  \end{CD}
\end{displaymath} 
whose vertical maps are isomorphisms. We consider a
$\mathfrak{S}$-linear isomorphism $\gamma: J \oplus L \rightarrow
J^{(\sigma)} \oplus L^{(\sigma)}$ which makes the following diagram
commutative
\begin{displaymath}
  \xymatrix{
    J \oplus L \ar[rr]^{\tau} \ar[rd]_{\gamma}& & M^{(\sigma)}\\
    &  J^{(\sigma)} \oplus L^{(\sigma)} .\ar[ru]_{\tau_1^{(\sigma)}} &
  }
\end{displaymath} 
The existence of $\gamma$ is implied by the following general
property. Let $N$ be a finitely generated module over a local ring
$A$. Let $F_1$ and $F_2$ be two free $A$-modules of the same rank
equipped with $A$-linear epimorphisms $\tau_1: F_1 \rightarrow N$ and
$\tau_2: F_2 \rightarrow N$. Then there exists an isomorphism
$\gamma_{12}: F_1 \rightarrow F_2$ such that we have $\tau_2 \circ
\gamma_{12} = \tau_1$.

We set $Q := J \oplus L$ and $\phi := \gamma \circ (E\id_J + \id_L): J
\oplus L \rightarrow J^{(\sigma)} \oplus L^{(\sigma)}$. Then the pair
$(Q, \phi)$ is a Breuil window relative to $\mathfrak{S}\to R$. We
have a commutative diagram
\begin{displaymath}
  \begin{CD}
    Q @>{\tau_1}>> M\\
    @V{\phi}VV  @V{\varphi}VV\\
    Q^{(\sigma)} @>{\tau_1^{(\sigma)}}>> M^{(\sigma)}.\\
  \end{CD}
\end{displaymath}
Hence $\tau_1$ is a surjection from $(Q,\phi)$ to $(M,\varphi)$. It is
obvious that the kernel of $\tau_1:(Q,\phi) \to (M,\varphi)$ is again
a Breuil module $(Q',\phi')$ relative to $\mathfrak{S}\to R$. We
obtain a short exact sequence:
\begin{displaymath}
  0 \rightarrow (Q',\phi') \rightarrow (Q,\phi) \rightarrow (M,\varphi)
  \rightarrow 0.
\end{displaymath}
Thus (ii) holds. 

Next we prove (iii).  We have seen above that for any Breuil module
$(M,\varphi)$ relative to $\mathfrak{S}\to R$ there is a Breuil window $(Q,\phi)$ relative to $\mathfrak{S}\to R$ and an epimorphism
$(Q,\phi) \rightarrow (M,\varphi)$. We remark that our argument uses
only the properties that $\varphi: M \rightarrow M^{(\sigma)}$ is
injective and that its cokernel is annihilated by $E$. 

In the
situation (iii) we choose a surjection $(\tilde{Q}, \tilde{\phi})
\rightarrow  (\tilde{M}, \tilde{\varphi})$ from a Breuil window $(\tilde{Q}, \tilde{\phi})$ relative to $\mathfrak{S}\to R$. 
We form the fibre product of $\mathfrak{S}$-modules $N = M \times_{\tilde{M}}
\tilde{Q}$. The functor which associates to an $\mathfrak{S}$-module
$L$ the $\mathfrak{S}$-module $L^{(\sigma)}$ is exact and therefore respects 
fibre products. We obtain a $\mathfrak{S}$-linear map $\psi : N \rightarrow
N^{(\sigma)}$ which is compatible with $\varphi$ and
$\tilde{\phi}$. Clearly $\psi$ is injective and its cokernel is annihilated by $E$. Therefore there is a
surjection $(Q, \phi) \rightarrow (N,\psi)$ from a Breuil window $(Q, \phi)$ relative to $\mathfrak{S}\to R$. We
deduce the existence of a commutative diagram 
\begin{displaymath}
   \begin{CD}
(Q,\phi) @>>> (M,\varphi)\\
@VVV  @VVV\\
(Q,\tilde{\phi}) @>>> (\tilde{M}, \tilde{\varphi}).
   \end{CD}
\end{displaymath}
As remarked above the kernels of the horizontal arrows are Breuil
modules relative to $\mathfrak{S}\to R$. This implies that (iii) holds.

To prove (iv) we consider the short exact sequence of (ii).  
As $\text{Coker}(\phi)$ and $\text{Coker}(\phi')$ are free $R$-modules
and as we have a short exact sequence $0\to \text{Coker}(\phi')\to
\text{Coker}(\phi)\to\text{Coker}(\varphi)\to 0$ of $R$-modules, we
get that (iv) holds as well.\endproof

\subsection{Proof of Theorem \ref{VZ2t}} 
We prove Theorem \ref{VZ2t}. Let $H$ be a finite flat, commutative
group scheme of $p$ power order over $R$. Due to a theorem of Raynaud
 (see \cite{BBM}, Theorem 3.1.1),
$H$ is the kernel of an isogeny of $p$-divisible groups over $R$
\begin{displaymath}
  0 \rightarrow H \rightarrow G^\prime \rightarrow G \rightarrow 0.
\end{displaymath}
Let $(Q',\phi')$ and $(Q, \phi)$ be the Breuil windows relative to
$\mathfrak{S}\rightarrow R$ which correspond by Theorem \ref{VZ1t} to
the $p$-divisible groups $G'$ and $G$. Let $(Q',\phi')\rightarrow
(Q,\phi)$ be the morphism that corresponds to the isogeny $G'\to
G$. This morphism is an isogeny i.e., it is a monomorphism and its
cokernel is annihilated by a power of $p$ (as $G'\rightarrow G$ is an
isogeny). Then it is immediate that the cokernel $(M,\varphi)$ of
$(Q',\phi')\rightarrow (Q,\phi)$ is a Breuil module relative to
$\mathfrak{S}\rightarrow R$. One can check that $(M,\phi)$ is
independent of the chosen resolution of $H$ and that the association
$H \mapsto (M,\varphi)$ is a functor.

Conversely let $(M, \phi)$ be a Breuil module relative to
$\mathfrak{S}\rightarrow R$. By Proposition \ref{P2} (ii) we have a
short exact sequence $0\to (Q',\phi')\to (Q,\phi)\to (M,\varphi)\to
0$, where $(Q',\phi')$ and $(Q,\phi)$ are Breuil windows relative to
$\mathfrak{S} \rightarrow R$. By Theorem \ref{VZ1t}, the monomorphism
$(Q',\phi')\to (Q,\phi)$ gives rise to an isogeny of $p$-divisible
groups $G' \rightarrow G$. Based on this and Proposition \ref{P2} (iii) we obtain a functor which associates to
$(M, \phi)$ the kernel of the isogeny $G^\prime\to G$. This is a
quasi-inverse to the functor of the previous paragraph. Thus Theorem
\ref{VZ2t} holds.\endproof

\smallskip

\hbox{Adrian Vasiu,\;\;\;Email: adrian@math.binghamton.edu}
\hbox{Address: Department of Mathematical Sciences, Binghamton University,}
\hbox{Binghamton, New York 13902-6000, U.S.A.}

\medskip
\hbox{Thomas Zink,\;\;\;Email: zink@math.uni-bielefeld.de}
\hbox{Address: Fakult\"at f\"ur Mathematik, Universit\"at Bielefeld,} \hbox{P.O. Box
100 131, D-33 501 Bielefeld, Germany.}

\end{document}